\documentclass[a4paper,12pt]{article}
\usepackage{amsmath}
\usepackage{amsthm}
\usepackage{amssymb}
\usepackage{amsfonts}
\usepackage{enumerate}
\usepackage[all]{xy}

\normalsize

\newtheoremstyle{theoremstyle}
  {10pt}      
  {5pt}       
  {\itshape}  
  {}          
  {\bfseries} 
  {:}         
  {.5em}      
  {}          

\newtheoremstyle{examplestyle}
  {10pt}      
  {5pt}       
  {}          
  {}          
  {\bfseries} 
  {:}         
  {.5em}      
  {}          

\theoremstyle{theoremstyle}
\newtheorem{theorem}{Theorem}[section]
\newtheorem{lemma}[theorem]{Lemma}
\newtheorem{proposition}[theorem]{Proposition}
\newtheorem*{proposition*}{Proposition}
\newtheorem{corollary}[theorem]{Corollary}
\newtheorem*{corollary*}{Corollary}

\theoremstyle{examplestyle}
\newtheorem{example}[theorem]{Example}
\newtheorem{definition}[theorem]{Definition}
\newtheorem{remark}[theorem]{Remark}
\newtheorem{remark*}{Remark}

\newcommand{\comment}[1]{}

\newcommand{\rays}{{\Delta(1)}}
\newcommand{\mrays}{$\rays$}
\newcommand{\weildivisors}{{\mathbb{Z}^\rays}}
\newcommand{\mweildivisors}{$\weildivisors$}
\newcommand{\sh}[1]{\mathcal{#1}}
\newcommand{\msh}[1]{$\sh{#1}$}
\newcommand{\spec}[1]{\operatorname{spec}(#1)}
\newcommand{\Spec}[1]{\operatorname{\mathbf{spec}}(#1)}
\newcommand{\dlim}{\underset{\rightarrow}{\lim}}
\newcommand{\mcoord}[1]{$\coord{#1}$}
\newcommand{\coord}[1]{k[{#1}_M]}
\newcommand{\morb}[1]{$\orb{#1}$}
\newcommand{\orb}[1]{\operatorname{orb}(#1)}
\newcommand{\fitt}[1]{\operatorname{Fitt}_2(#1)}

\title{Graded Rings and Equivariant Sheaves on Toric Varieties}
\author{Markus Perling\\ \footnotesize Department of Mathematics\\[-9pt]
        \footnotesize University of Kaiserslautern, Germany\\[-9pt] 
        \footnotesize \it email: \tt perling@mathematik.uni-kl.de}
\date{October 2001}

\begin{document}

\maketitle

\begin{abstract}
In this note we derive a formalism for describing equivariant sheaves over toric
varieties. This formalism is a generalization of a correspondence due to Klyachko,
which states that equivariant vector bundles on toric varieties are equivalent to
certain sets of filtrations of vector spaces. We systematically construct the theory
from the point of view of graded ring theory and this way we considerably clarify
earlier constructions of Kaneyama and Klyachko.
 We also connect the formalism to the theory of fine-graded
modules over Cox' homogeneous coordinate ring of a toric variety. As an application we
construct minimal resolutions of equivariant vector bundles of rank two on toric surfaces.
\end{abstract}

\section{Introduction}

In \cite{Kly90} Klyachko showed that equivariant vector bundles on a toric variety
are equivalent to certain sets of filtrations of vector spaces. In the preprint \cite{Kly91}
this equivalence was extended to equivariant torsion free sheaves on smooth toric varieties
and
applied to moduli problems for bundles over $\mathbb{P}_2$. It was claimed there that
results similar to that of \cite{Kly90} concerning vector bundles, are valid
for torsion free sheaves, but unfortunately no proofs were given. The aim of this note
is to pick up some of the ideas of \cite{Kly91} and to deliver proper definitions of the
filtration-formalism for equivariant torsion free sheaves and to generalize it to
arbitrary equivariant coherent sheaves on not necessarily smooth toric varieties.

One basic idea is that if one considers an equivariant quasicoherent sheaf \msh{E}
over the toric variety $X_\Delta$ associated to some fan $\Delta$ (see Section
\ref{basicfacts} for notation), then the action of the torus $T$ on the affine
$T$-invariant open subsets $U_\sigma$, $\sigma
\in \Delta$, induces isotypical decompositions of the modules of sections of \msh{E} into
$T$-eigenspaces:
\begin{equation*}
\Gamma(U_\sigma, \sh{E}) = \bigoplus_{m \in M}\Gamma(U_\sigma, \sh{E})_m,
\end{equation*}
where $M$ denotes the character group of the torus which we
identify with $\mathbb{Z}^{\dim T}$.
If \msh{E} is coherent, the eigenspaces $\Gamma(U_\sigma, \sh{E})_m$ are
finite-dimensional vector spaces. For each $m, m' \in M$ the module structure
over the coordinate ring \mcoord{\sigma} of $U_\sigma$ induces maps
\begin{align*}
\Gamma(U_\sigma, \sh{E})_m & \longrightarrow \Gamma(U_\sigma, \sh{E})_{m'} \\
e & \mapsto \chi(m' - m) \cdot e
\end{align*}
by multiplication with the character $\chi(m' - m)$.
This multiplication map exists if and only if $m' - m$ is an element of the semigroup
$\sigma_M = \check{\sigma} \cap M$.
So the semigroup $\sigma_M$ induces in a natural way a preorder on $M$ by setting
$m \leq m'$ iff $m' - m \in \sigma_M$. The set of
vector spaces $\Gamma(U_\sigma, \sh{E})_m$ and characters $\chi(m)$ forms a directed family
of vector spaces with respect to this preorder, and we will denote such data a
{\it $\sigma$-family}.
We obtain for each $\sigma \in \Delta$ such a family, and vice versa, given a set of
$\sigma$-families for all $\sigma \in \Delta$, we obtain a system of sheaves
$\sh{E}_\sigma$ over each $U_\sigma$ which, if certain compatibility conditions are
fulfilled, glue to a global sheaf over the toric variety $X_\Delta$.
Such a compatible set of $\sigma$-families will be called a {\it $\Delta$-family}.
Ultimately, we will arrive at an equivalence of categories between equivariant
quasicoherent sheaves over $X_\Delta$ and $\Delta$-families.
If \msh{E} is torsion free, all the maps in its associated $\Delta$-family
are injective and one can formulate the result in terms of multifiltrations of a
certain limit vector space $\mathbf{E}^0$. This is essentially the formulation which
Klyachko gave in \cite{Kly91} and which generalized the earlier filtrations for
equivariant vector bundles of \cite{Kly90}. We show that the constructions of
\cite{Kly90} and \cite{Kly91} simplify considerably if they are derived systematically
from the point of view of graded ring theory. This point of view also has been adopted
earlier by Kaneyama in \cite{Kan1, Kan2} for classification of equivariant vector
bundles on smooth toric varieties and especially on $\mathbb{P}_n$.

Moreover, we will enhance these constructions by showing how fine graded modules over the
homogeneous coordinate ring of a toric variety (\cite{Cox}) can be incorporated.
As an application we construct minimal resolutions of equivariant vector bundles of
rank 2 on smooth toric surfaces in terms of equivariant Euler type short exact sequences.
In a forthcoming paper we will extend this application
to a full classification of such vector bundles.

{\bf Plan of the paper:}
In Section \ref{basicfacts}, we will recall some basic notions concerning toric varieties.
Section \ref{gradings} collects general facts of graded ring theory for reference.
In Section \ref{homcoordring} we briefly desribe the construction
of homogeneous coordinate rings over toric varieties due to Cox. The main part of this
note is Section \ref{equivariantsheaves}, where we develop our formalism for equivariant
sheaves, and we conclude with an application in Section \ref{resolutions}.

{\bf Acknowledgements:} I would like to thank Prof. G. Trautmann and Dr. J. Zintl for many
discussions on the subject.

\section{Basic Facts for Toric Varieties}
\label{basicfacts}

By an {\it algebraic variety over an algebraically closed field $k$}, we understand a
separated scheme of finite type over $\Spec{k}$, which we assume to be reduced,
if not stated otherwise.
A toric variety $X$ is a normal variety which contains an algebraic torus
$T$ as an open dense subset such that the torus multiplication extends to
an action of the algebraic group $T$ on $X$. For background on toric varieties we refer
to standard literature such as \cite{MiyakeOda}, \cite{Oda} and \cite{Fulton}. We recall
from \cite{MiyakeOda}, \S\S\ 5.3, 5.4 and \cite{Fulton} some basic facts which in
the sequel will frequently be used. A toric variety
$X$ is defined by a fan $\Delta$ contained in the real vector space $N_\mathbb{R} \cong N
\otimes_{\mathbb{Z}} \mathbb{R}$ of a lattice $N \cong \mathbb{Z}^n$ and is denoted by
$X = X_\Delta$.
Let $M$ be the lattice dual to $N$ and let $\langle \, , \rangle : M \times N \rightarrow
\mathbb{Z}$ be the natural pairing. This pairing extends in a natural way to the scalar
extensions $M_\mathbb{R} := M \otimes_\mathbb{Z} \mathbb{R}$
and $N_\mathbb{R}$. Elements of $M$ are denoted by $m$, $m'$, etc. if written additively,
and by $\chi(m)$, $\chi(m')$, etc. if written multiplicatively, i.e. $\chi(m + m') = \chi(m)
\chi(m')$. The lattice $M$ is the natural group of characters of the torus
$T = \operatorname{Hom}(M,k^*) \cong (k^*)^n$.

A cone $\sigma$ of the fan $\Delta$ is a convex rational polyhedral cone contained in
$N_\mathbb{R}$. For these cones the following standard notation will be used:

\begin{itemize}
\setlength{\itemsep}{-5pt}
\item cones are denoted by small Greek letters $\rho$, $\sigma$, $\tau$, etc.,
the natural order among cones is denoted by $\tau < \sigma$,
\item $\Delta(i) := \{\sigma \in \Delta \mid \dim \sigma = i\}$ the set of
all cones of fixed dimension $i$,\\
$\sigma(i) := \{ \tau \in \Delta(i) \mid \tau < \sigma\}$, elements of $\Delta(1)$ are called {\it rays},
\item $\vert \Delta \vert := \bigcup_{\sigma \in \Delta} \sigma$ is the {\it support} of the fan $\Delta$,
\item $n(\rho)$ denotes the primitive lattice element spanning the ray $\rho$,
\item $\check{\sigma} := \{m \in M_\mathbb{R} \mid \langle m, n \rangle \geq 0 \text{ for all
$n \in \sigma$}\} = \bigcap_{\rho \in \sigma(1)} H_{n(\rho)}$ is the cone {\it dual}
to $\sigma$; $H_{n(\rho)}$ is the positive half space defined by the primitive lattice vector $n(\rho)$
and is the same as $\check{\rho}$,
\item $\sigma^\bot = \{m \in M_\mathbb{R} \mid m(n) = 0 \text{ for all } n \in \sigma \}$,
\item $\sigma_M := \check{\sigma} \cap M$ is the subsemigroup of $M$ associated to $\sigma$.
$\sigma_M$ is a
finitely generated semigroup which is saturated (i.e. if $c \cdot m \in \sigma_M$ with $0 < c \in \mathbb{N}$,
and $m \in M$, then $m \in \sigma_M$) and which generates $M$ as a group; $\sigma_M$ is
often denoted by $S_\sigma$ (cf. \cite{Fulton}).\\
$\sigma_M^\bot := \sigma^\bot \cap M$ is the unique maximal subgroup of $M$ contained in $\sigma_M$.
\end{itemize}
Dualizing a cone $\sigma$ induces an order-reversing correspondence between the faces of
$\sigma$ and the faces of $\check{\sigma}$. This correspondence is given by $\tau
\longleftrightarrow \check{\sigma} \cap \tau^\bot$.
A cone $\sigma$ of $\Delta$ corresponds to an open affine subset $U_\sigma \subset X_\Delta$ which
is an affine toric variety. This correspondence is natural in the sense that $U_\tau \subset
U_\sigma$ for $\tau < \sigma$ and $U_\sigma \cap U_{\sigma'} = U_{\sigma \cap \sigma'}$
for $\sigma, \sigma' \in \Delta$.
The coordinate ring $k[U_\sigma]$ of $U_\sigma$ is the semigroup ring
\begin{equation*}
\coord{\sigma} = \bigoplus_{m \in \sigma_M} k \cdot \chi(m).
\end{equation*}
The torus $T$ acts on \mcoord{\sigma} via its dual
action by $t . \sum a_i \chi(m_i) = \sum a_i t^{-m_i}\chi(m_i)$.

\paragraph{$T$-invariant subvarieties of $U_\sigma$:}

The ring \mcoord{\sigma} has the structure of a $\sigma_M$-graded (and hence $M$-graded) ring
and there exists a natural bijection  between $\sigma_M$-homogeneous ideals of \mcoord{\sigma}
and semigroup ideals of $\sigma_M$ as follows. A $\sigma_M$-homogeneous ideal $\mathfrak{a}$
defines the semigroup ideal $I = \{m \in \sigma_M \mid \chi(m) \in \mathfrak{a}\}$ and the ideal
$\mathfrak{a}$ is determined by $I$ as
\begin{equation*}
\mathfrak{a} = \mathfrak{a}_I =  \bigoplus_{m \in I} k \cdot \chi(m).
\end{equation*}
Then the closed (not necessarily reduced) subvariety $\mathbf{V}(\mathfrak{a})
\subset U_\sigma$ has the coordinate ring
\begin{equation*}
\coord{\sigma}/\mathfrak{a}_I = \bigoplus_{m \in \sigma_M \setminus I}k \cdot \chi(m).
\end{equation*}
The ideal $\mathfrak{a}_I$ is a prime ideal iff there is a unique face $\tau$ of $\sigma$
such that $\sigma_M \setminus I = \sigma_M \cap \tau^\bot$. On the other hand, a face $\tau$
of $\sigma$ determines a $\sigma_M$-homogeneous prime ideal $\mathfrak{p}_\tau$ by this
formula. This establishes that the faces of $\sigma$ correspond 1-1 both with the
$\sigma_M$-homogeneous prime ideals of \mcoord{\sigma} and with the $T$-invariant integral
subvarieties of $U_\sigma$. We denote by $V_\sigma(\tau) \subset U_\sigma$ the subvariety
determined by the face $\tau$. Its coordinate ring is
\begin{equation*}
\bigoplus_{m \in \sigma_M \cap \tau^\bot} k \cdot \chi(m).
\end{equation*}
There is a unique largest prime ideal $\mathfrak{p}_\sigma$ of \mcoord{\sigma} which is associated
to the semigroup ideal generated by the complement of the maximal subgroup $\sigma_M^\bot$ in
$\sigma_M$:
\begin{equation*}
\mathfrak{p}_\sigma = \bigoplus_{\sigma_M \setminus \sigma_M^\bot} k \cdot \chi(m).
\end{equation*}
The corresponding variety is the unique minimal (and hence the unique closed) orbit of
the torus action on $U_\sigma$:
\begin{equation*}
\orb{\sigma} := V_\sigma(\sigma) = \spec{k[\sigma^\bot_M]}.
\end{equation*}
\morb{\sigma} is a torus whose dimension is equal to the codimension of $\sigma$
in $N_\mathbb{R}$. The stabilizer group of \morb{\sigma} is the torus
associated to the character group $M / \sigma_M^\bot$.

Let $\tau < \sigma$, then there exists a lattice element $m_\tau$ in the interior of $\tau^\bot
\cap \check{\sigma}$ such that $\tau_M = \sigma_M + \mathbb{Z}_{\geq 0} (-m_\tau)$ and
$\tau_M^\bot = \tau^\bot_M \cap \sigma_M + \mathbb{Z}_{\geq 0}(-m_\tau)$;
this implies that
\begin{equation*}
\coord{\tau}=\coord{\sigma}_{\chi(m_\tau)} \text{ \quad and \quad } k[\tau^\bot_M] =
k[\tau_M^\bot \cap \sigma_M]_{\chi(m_\tau)} = k[V_\sigma(\tau)]_{\chi(m_\tau)}.
\end{equation*}
where the index denotes localization.
Thus for each $\tau < \sigma$
\begin{itemize}
\setlength{\itemsep}{-5pt}
\item the orbit \morb{\tau} can be identified with the intersection $V_\sigma(\tau) \cap U_\tau$,
\item $V_\sigma(\tau)$ is precisely the closure of \morb{\tau} in $U_\sigma$, and
\item $\orb{\sigma}$ is contained in $V_\sigma(\tau)$.
\end{itemize}
For any cone $\tau \in \Delta$ we can consider the closure $V(\tau)$ of $\orb{\tau}$ in
$X_\Delta$ and for any $\sigma > \tau$ we identify $V(\tau) \cap U_\sigma$ with $V_\sigma(\tau)$.
There is an orbit decomposition:
\begin{equation*}
V(\tau) = \bigcup_{\sigma > \tau} V_\sigma(\tau) = \coprod_{\sigma > \tau} \orb{\sigma}.
\end{equation*}

If $\rho \in \rays$, its associated closed subvariety $V(\rho)$ has codimension 1 in $X_\Delta$
and is an irreducible, $T$-invariant divisor of $X_\Delta$, denoted $D_\rho$. The group
$\operatorname{TWDiv}(X_\Delta)$ of $T$-invariant Weil divisors of $X_\Delta$ is the group freely
generated over $\mathbb{Z}$ by the divisors $D\_rho$, $\rho \in \rays$.
$\operatorname{TWDiv}(X_\Delta)$ thus is often identified with $\mathbb{Z}^\rays$
which we will use in the sequel as standard notation.
\begin{itemize}
\item By $A$ we denote the Chow group $A_{n - 1}(X_\Delta)$.
\end{itemize}
By \cite{MiyakeOda}, Proposition 6.1 and \cite{Fulton}, \S 3.4, there exists a short
exact sequence:
\begin{equation}
\label{maindiagram}
0 \longrightarrow M \longrightarrow \weildivisors \longrightarrow A \longrightarrow 0,
\end{equation}
where the embedding $0 \longrightarrow M \longrightarrow \weildivisors$ is defined by
$m \mapsto \langle m, n(\rho) \rangle$. The group $A$ is thus generated by
the $T$-invariant Weil divisors.
This sequence exists for all toric varieties whose fan is not contained in a proper
subspace of $N_\mathbb{R}$. In the sequel we will always assume the existence of the
sequence.

\section{Preliminaries on Graded Rings}
\label{gradings}

In this section we collect some standard facts on graded commutative rings and modules.
Let $G$ be an abelian group and let $R$ be a {\it $G$-graded} ring.
A module $M$ over a $G$-graded ring $R$ is called a {\it graded $R$-module} if $M$ has a
decomposition into a direct sum $M = \bigoplus_{g \in G}M_g$ of abelian subgroups $M_g$
of $M$ such that $r \cdot m \in M_{gg'}$ for all $r \in R_g, m \in M_{g'}$.
If $M$ is an ideal of $R$, it is also called {\it homogeneous}.
We denote by $M(g)$ the degree shift such that $M(g)_{g'} = M_{gg'}$.
Let $M, N$ be two $G$-graded $R$-modules.
A homomorphism $f : M \longrightarrow N$ of two $G$-graded $R$-modules $M$ and $N$ is
a {\it graded} or {\it homogeneous homomorphism} if $f(M_{g}) \subset N_{g}$ for
all $g \in G$.
We can define a graded structure on the tensor product $M \otimes_R N$ as follows:
denote by $T := M  \otimes_\mathbb{Z} N$ the tensor product of $M$ and $N$ over $\mathbb{Z}$.
If we consider $\mathbb{Z}$ as a trivially $G$-graded ring (i.e. $\mathbb{Z}_1 = \mathbb{Z}$
and $\mathbb{Z}_g = 0$ for all $1 \neq g \in G$) then $T$ has a $G$-grading by
setting $T_g$ to be the additive group generated by the set $\{ m \otimes n \mid m \in M_{g'},
n \in N_{g''} \text{ such that } g'g'' = g \}$ for all $g \in G$.
Further let $K$ be the $G$-graded submodule of $T$ generated by $\{ (rm) \otimes n - m \otimes
(rn) \mid m \in M, n \in N, r \in R \}$.
There is an isomorphism $T / K \cong M \otimes_R N$
such that $M \otimes_R N$ acquires the structure of a $G$-graded module.

\begin{remark}
If $R$ is trivially $G$-graded then $M \otimes_R N$ can be endowed directly with a
graded structure by setting $(M \otimes_R N)_g$ to be the additive group generated by
$\{ m \otimes_R n \mid m \in M_{g'}, n \in N_{g''} \text{ such that } g'g'' = g \}$.
\end{remark}

A {\it morphism of graded rings} is a pair of morphisms $(\psi, \chi) : (R, G)
\longrightarrow (R', G')$ such that $\psi(R_g) \subset R'_{\chi(g)}$.
Such a morphism $(\psi, \chi)$
induces a $G'$-grading on $R$ by setting $R_{g'} := \bigoplus_{g \in \chi^{-1}(g')}
R_g$. Also for a given graded $G$-module $M$ the grading can also be changed
to a $G'$-grading. Thus:

\begin{definition}
\label{gradscalext}
Let $(\psi, \chi) : (R, G) \longrightarrow (R', G')$ be a morphism of commutative graded
rings and let $M$ be a $G$-graded $R$-module. Then we call the $G'$-graded $R'$-module $M
\otimes_R R'$ the {\it graded scalar extension}.
\end{definition}

\begin{definition}
Let $R$ be a $G$-graded ring.
A {\it maximal homogeneous ideal} is a homogeneous
ideal of $R$ which is not contained in any other homogeneous ideal of $R$.
The intersection of all maximal homogeneous ideals is called the
{\it homogeneous Jacobson radical} $\mathfrak{R}_G$.
\end{definition}

Note that a maximal homogeneous ideal is not necessarily a maximal ideal.

\begin{theorem}[Graded Version of Krull's Theorem]
\label{gradedkrull}
Let $R$ be a graded ring. Then each proper homogeneous ideal is contained in a
maximal homogeneous ideal.
\end{theorem}

\begin{proof}
Clearly the set of homogeneous ideals is partially ordered by inclusion.
Let $\{\mathfrak{g}_i\}_{i \in I}$ be a chain of homogeneous ideals. Then
the ideal $\bigcup_{i \in I}\mathfrak{g}_i$ is also homogeneous.
So the homogeneous ideals are even {\it inductively} ordered, hence we can apply
Zorn's lemma and obtain that there exist maximal homogeneous ideals.
\end{proof}

As a consequence we have the

\begin{corollary}
\label{maxunit}
\label{maxnullfree}
\begin{enumerate}[(i)]
\item Every homogeneous non-unit of $R$ is contained in a maximal homogeneous ideal.
\item If the zero ideal is the unique maximal homogeneous ideal of a $G$-graded ring $R$,
then each finitely generated graded module over $R$ is free.
\item
Let $R$ be a $G$-graded ring and assume that there exists
a unique maximal homogeneous ideal $\mathfrak{m}$ in $R$. The ring $R / \mathfrak{m}$
is then also graded and each finitely generated graded module over $R / \mathfrak{m}$
is free.
\end{enumerate}
\end{corollary}

\begin{theorem}[Graded Nakayama's lemma]
\label{gradnak}
Let $R$ be a $G$-graded ring, and let $M$ be a finitely
generated graded $R$-module. Let $\mathfrak{a} \subset \mathfrak{R}_G$ be a graded
ideal such that $\mathfrak{a}M = M$. Then $M = 0$.
\end{theorem}

\begin{proof}
Assume that $M \neq 0$, and let $u_1, \dots, u_n$ be a minimal set of generators
of $M$. Then $u_n \in \mathfrak{a}M$, hence there exists an equation of the form
$u_n = \sum_{i = 1}^n a_i u_i$ where $a_i \in \mathfrak{a}$. Write
$(1 - a_n)u_n = \sum_{i = 1}^{n - 1} a_i u_i$. Since $a_n \in \mathfrak{a}$,
by Corollary \ref{maxunit}, $1 - a_n$ is a unit in $R$. Hence $u_n \in \operatorname{span}_A
\{u_i\}_{i = 1, \dots, n - 1}$ which contradicts the minimality of the generating set.
\end{proof}

\begin{corollary}
\label{locgradprojfree}
Assume that there exists
a unique maximal homogeneous ideal $\mathfrak{m}$ in $R$. Let $M$ be a finitely generated
graded $R$-module. If $M$ is projective, then $M$ is free.
\end{corollary}

\section{The Quotient Construction of Cox}
\label{homcoordring}
\label{coxqcoh}

We briefly recall a construction of Cox and others (see \cite{Cox} and references
therein) who, generalizing the notion of the affine cone over a projective space,
gave a construction of toric varieties as categorical quotients of quasiaffine
varieties by diagonalizable groups.
The starting point is sequence
(\ref{maindiagram}), from which we obtain, by applying
$\operatorname{Hom}_\mathbb{Z}(\_, k^*)$:
\begin{equation*}
1 \longrightarrow G \longrightarrow \hat{T} \longrightarrow T \longrightarrow 0,
\end{equation*}
where $\hat{T} := (k^*)^{\rays}$ is the maximal torus contained in $k^\rays$ and
$G = \operatorname{Hom}_\mathbb{Z}(A, k^*)$ is a diagonalizable group which is
not necessarily irreducible.

Let $C$ be the cone which is generated over
$\mathbb{R}_{\geq 0}$ by the standard basis $\{e_\rho\}_{\rho \in \rays}$ of
$\mathbb{R}^{\rays}$. The associated toric variety then can be identified with the affine
space $k^{\rays}$, together with the action of the embedded torus $\hat{T}$.
We define $\hat{\Delta}$ to be the subfan of the fan generated by $C$ such that
the primitive elements $e_{\rho_1}, \dots, e_{\rho_l}$
span a cone $\hat{\sigma}$ in $\hat{\Delta}$ if and only if $n(\rho_1), \dots n(\rho_l)$
span a cone $\sigma$ in $\Delta$. Because $\hat{\Delta}$ is a subfan of $C$, the
corresponding toric variety $\hat{X}$ is an open subset of $k^\rays$ and thus
quasi-affine.

Denote $S := k[x_\rho \mid \rho \in \rays]$ the coordinate ring of $k^{\rays}$. An
affine $k^\rays$-invariant open subset $U_{\hat{\sigma}}$ of $\hat{X}$ is the complement
in $k^\rays$ of the union of all $T$-invariant divisors $D_{\hat{\rho}}$ such that
$\hat{\rho}$ is not in $\hat{\sigma}(1)$. The defining equation for such a union
of hypersurfaces is given by the monomial
\begin{equation*}
x^{\hat{\sigma}} := \prod_{\rho \in \rays \setminus \sigma(1)} x_\rho,
\end{equation*}
and the coordinate ring of $U_{\hat{\sigma}}$ is the localization
$S_{x^{\hat{\sigma}}}$. So the complement of $\hat{X}$ in $k^{\rays}$ is given
by the ideal
\begin{equation*}
B := \langle x^{\hat{\sigma}} \mid \hat{\sigma} \in \hat{\Delta} \rangle.
\end{equation*}
and has {\it at least codimension two}. $B$ is called the {\it irrelevant ideal}.

By dualizing the inclusion of $M$ into \mweildivisors\ in diagram (\ref{maindiagram}),
we obtain a linear map from \mweildivisors\ to $N$ which maps $e_\rho$ to $n(\rho)$ for each
$\rho \in \rays$.
By construction, this map induces a map of fans $\hat{\Delta} \longrightarrow \Delta$
and thus an equivariant map $\hat{X} \overset{\pi}{\longrightarrow} X$ of toric
varieties. The morphism $\pi$ maps $U_{\hat{\sigma}}$ onto $U_\sigma$ for all
$\sigma \in \Delta$ and has a natural interpretation as a categorical quotient of $\hat{X}$
by $G$. Namely, for each $\hat{T}$-invariant affine open subset $U_{\hat{\sigma}}$ of
$\hat{X}$ the dual action of $G$ on $S_{x^{\hat{\sigma}}}$ induces a grading by the
character group $A$ of $G$:
\begin{equation*}
S_{x^{\hat{\sigma}}} = \bigoplus_{\alpha \in A} (S_{x^{\hat{\sigma}}})_\alpha.
\end{equation*}
The ring $S_{x^{\hat{\sigma}}}$ carries also a natural grading by the group \mweildivisors,
and the $A$-grading is given by the surjective group homomorphism $\weildivisors
\twoheadrightarrow A$.
The $G$-invariant regular functions over $U_{\hat{\sigma}}$ then form the
subring of degree-zero elements, $S_{(x^{\hat{\sigma}})}$.
It was shown explicitly in \cite{Cox} that $S_{(x^{\hat{\sigma}})}$ is isomorphic to
$k[\sigma_M]$ which implies that $\pi$ is a universal categorical quotient of
$\hat{X}$ by $G$ in the sense of \cite{GIT}.
The variables $x_\rho$, $\rho \in \rays$ of $S$ serve, modulo the action of $G$,
as global coordinates for $X$ and, together with the $A$-grading
$S = \bigoplus_{\alpha \in A} S_\alpha$, the ring $S$ is called the
{\it homogeneous coordinate ring} of $X$.

\paragraph{Quasi-coherent sheaves:}

Given an $A$-graded $S$-module $F$, then for each $\sigma \in \Delta$ one considers the
$S_{x^{\hat{\sigma}}}$-module $F_{x^{\hat{\sigma}}}$.
Taking degree zero yields an $S_{(x^{\hat{\sigma}})}$-module $F_{(x^{\hat{\sigma}})}$
and thus applying the Grothendieck functor we obtain a quasicoherent sheaf
$\tilde{F}_{(x^{\hat{\sigma}})}$ on $U_\sigma$. These
sheaves glue together to a sheaf on $X$ which is denoted $\tilde{F}$. This
sheafification operation defines an exact functor from the category of graded $S$-modules
to the category of quasicoherent sheaves over $X$.

An important class of sheaves is coming from free $S$-modules of rank 1,
which are of the form $S(\alpha)$ for some $\alpha \in A$, where $S(\alpha)$ denotes
the {\it degree shift} of $S$ by $\alpha$, i.e. $S(\alpha)_\beta = S_{\alpha + \beta}$.
If $D = \sum n_Y Y$
is a Weil divisor representing $\alpha$, a sheaf
 $\sh{O}_X(D)$ can be defined
by: $U \mapsto \Gamma(U,\sh{O}_X(D)) := \{f \in k(X) \mid \upsilon_Y(f) \geq -n_Y \;
\; \forall \; Y \in U\}$ where $\upsilon_Y$ is the valuation of the generic
point of $Y$.
If $D$ is a Cartier divisor, the sheaf $\sh{O}_X(D)$ is the invertible sheaf associated to $D$
in the usual way, but in general $\sh{O}_X(D)$ is only a reflexive sheaf of rank 1.
By \cite{Reid1}, Theorem (3) there exists a bijection between
$A$ and the set of isomorphism classes of reflexive sheaves of rank 1.
This correspondence is conveniently compatible with the ring structure of $S$:
if $D$ is a Weil divisor representing $\alpha \in A$, then there is an isomorphism
$\widetilde{S(\alpha)}$ $\cong \sh{O}_X(D)$ and the global sections $\Gamma(X, \sh{O}_X(D))$
can be identified with $S_\alpha$, the homogeneous elements of degree zero of the module
$S(\alpha)$.
 The natural map
$\Gamma(X, \sh{O}_X(D)) \otimes_{\Gamma(X, \sh{O}_X)} \Gamma(X, \sh{O}_X(D'))
\longrightarrow \Gamma(X, \sh{O}_X(D + D'))$ then can naturally be identified with the ring
multiplication of $S$ (cf. \cite{Cox}, Prop. 1.1).

Let us denote by $\sh{O}_X(\alpha) := \widetilde{S(\alpha)}$ the distinguished
representative of the isomorphism class of sheaves $\sh{O}_X(D)$ where $D$ is a
representative of $\alpha$.
Then, in the other direction, for a quasicoherent sheaf \msh{F} on $X$ one can define a
graded $S$-module
\begin{equation*}
\Gamma_* \sh{F} := \bigoplus_{\alpha \in A} \Gamma\big(X, \sh{F} \otimes_{\sh{O}_X}
\sh{O}_X(\alpha)\big ).
\end{equation*}
$\Gamma_*$ is a functor from the category of quasi-coherent sheaves over $X$ to the
category of graded $S$-modules. In the following theorem we summarize some results
relating $A$-graded modules and quasicoherent sheaves.

\begin{theorem}[\cite{Cox}, \cite{Mustata1}]
\begin{enumerate}[(i)]
\item The map $F \mapsto \tilde{F}$ is a covariant additive exact functor
from the category of $A$-graded $S$-modules to the category of quasi-coherent $\sh{O}_X$-modules.

\item If $F$ is a finitely generated $A$-graded $S$-module, then $\tilde{F}$ is coherent.

\item\label{mustext} There exists a natural isomorphism
$(\Gamma_*\sh{F})\tilde{\ } \cong \sh{F}$. In particular, every quasi-coherent $\sh{O}_X$--module
\msh{F} is of the form $\tilde{F}$ for some $A$-graded $S$-module $F$. If \msh{F} is coherent,
then there exists
a finitely generated $A$-graded $S$-module $F$ with $\tilde{F} \cong \sh{F}$.
\end{enumerate}
\end{theorem}

(\ref{mustext}) has been proven by Cox for the case of simplicial toric varieties and by
Mustata (\cite{Mustata1}) in general.
We complement these results by the following remark:

\begin{proposition}
Let \msh{F} be a quasi-coherent sheaf on $X$. Then $\Gamma_*{\sh{F}} \cong
\Gamma(\hat{X}, \pi^*\sh{F})$.
\end{proposition}

\begin{proof}
By restricting the projection $\pi$ to any $U_{\hat{\sigma}}$ we get a map of affine
varieties $\pi\vert_{U_{\hat{\sigma}}} : U_{\hat{\sigma}} \longrightarrow U_\sigma$.
The induced map $\pi' : \coord{\sigma} \longrightarrow
k[U_{\hat{\sigma}}]$ yields a homomorphism of graded rings $(\pi', 0) : (\coord{\sigma}, M)
\longrightarrow (k[U_{\hat{\sigma}}], A)$, where $0 : M \longrightarrow A$ is just the zero
map. Let us denote $\Gamma(U_\sigma, \sh{F})$ by $F_\sigma$, then $\Gamma(U_{\hat{\sigma}},
\pi^* \sh{F})$ is isomorphic to $k[U_{\hat{\sigma}}] \otimes_{\coord{\sigma}} F_\sigma$,
where $k[U_{\hat{\sigma}}] \otimes_{\coord{\sigma}} F_\sigma$ can be considered as a graded
scalar extension (see \ref{gradscalext}) with $F_\sigma$ and \mcoord{\sigma} trivially
$A$-graded.
Thus the pullback decomposes into
homogeneous components:
\begin{align*}
\Gamma(U_{\hat{\sigma}}, \pi^* \sh{F}) & \cong \bigg (\bigoplus_{\alpha \in A}
(S_{x^{\hat{\sigma}}})_\alpha \bigg) \otimes_{\coord{\sigma}} F_\sigma \\
& \cong  \bigoplus_{\alpha \in A} \big ((S_{x^{\hat{\sigma}}})_\alpha
\otimes_{\coord{\sigma}} F_\sigma \big)\\
& \cong \bigoplus_{\alpha \in A} \Gamma\big (U_\sigma, \sh{O}_X(\alpha)\big )
\otimes_{\coord{\sigma}} \Gamma(U_\sigma, \sh{F}) \\
&  \cong \bigoplus_{\alpha \in A} \Gamma\big (U_\sigma, \sh{F} \otimes \sh{O}_X(\alpha)\big).
\end{align*}
Gluing then yields the result.
\end{proof}

\begin{corollary}
\label{pullbacktilde}
Consider the map $\hat{X} \overset{\pi}{\longrightarrow} X$. Then
$\big(\Gamma(\hat{X}, \pi^*\sh{F})\big)\tilde{\ } \cong \sh{F}$.
\end{corollary}

\section{Equivariant Sheaves}
\label{equivariantsheaves}

Let $G$ be an algebraic group which acts on a variety $X$. Following \cite{GIT}, \S 3, we call
a sheaf \msh{E} on $X$
{\it equivariant} or {\it $G$-linearized} if for each $g \in G$ considered as automorphism of
$X$, there is an isomorphism $\Phi_g : g^* \sh{E} \cong \sh{E}$ such that the following diagram
commutes for all $g, g' \in G$:

$$\xymatrix{
g^*{g'}^* \sh{E} \ar[rrd]_{g^*\Phi_{g'}} & & (g'g)^* \sh{E} \ar[ll]_\cong \ar@{-->}[d] \ar[rr]^{\Phi_{g'g}} & &  \sh{E} \\
& & {g}^* \sh{E} \ar[urr]_{\Phi_g}
}$$
Over a $G$-invariant open subset $U$ of $X$ there is an induced action of $G$ on
$\Gamma(U_, \sh{E})$. For $s \in \Gamma(U_\sigma, \sh{E})$ let $g^*s$ denote the
canonically lifted section of $g^*\sh{E}$
and let the dual action of $G$ on $\Gamma(U, \sh{E})$ be defined by
$g.s = \phi_{g^{-1}}(g^{-1 *}s)$.
This action is a linear representation of $G$ and is called the {\it dual action} of $G$ on
$\Gamma(U_, \sh{E})$.
We consider the case where $G$ is an algebraic torus over an algebraically closed field $k$.

\subsection{The Category of $\sigma$-Families}

Let \msh{E} be a quasi-coherent equivariant sheaf on an affine toric variety $U_\sigma$
and denote by $E^\sigma$ the \mcoord{\sigma}-module $\Gamma(U_\sigma, \sh{E})$. We consider
its isotypical decomposition:
\begin{equation*}
E^\sigma = \bigoplus_{m \in M} E^\sigma_m.
\end{equation*}
The $M$-graded module structure of $E^\sigma$ is determined by the mappings given
by multiplication with characters in $\sigma_M$:
\begin{align*}
\chi(m): E^\sigma_{m'} & \longrightarrow E^\sigma_{m + m'}\\
e & \mapsto \chi(m) \cdot e
\end{align*}
for all $m \in \sigma_M$ and all $m' \in M$. Thus there exists such a map between the
eigenspaces $E^\sigma_m$ and $E^\sigma_{m'}$ if and only if $m' - m \in \sigma_M$.
Let $F^\sigma$
be another $M$-graded \mcoord{\sigma}-module. Then a homogeneous morphism $\phi$ of
degree 0 between $E^\sigma$ and $F^\sigma$ decomposes into a sum
\begin{equation*}
\phi = \bigoplus_{m \in M} \phi_m
\end{equation*}
of $k$-vector space homomorphisms $\phi_m : E^\sigma_m \longrightarrow F^\sigma_m$
satisfying the compatibility condition
\begin{equation*}
\phi_{m + m'} \circ \chi(m) =  \chi(m) \circ \phi_{m'}
\end{equation*}
for all $m \in \sigma_M$ and all $m' \in M$. We will make these properties precise:

\begin{definition}
\label{preorderdef}
Let $U_\sigma$ be an affine toric variety. Then we define a relation $\leq_\sigma$ on $M$
by setting $m \leq_\sigma m'$ iff $m' - m \in \sigma_M$.
We write $m <_\sigma m'$ if $m \leq_\sigma m'$ but not $m' \leq_\sigma m$.
\end{definition}

It is then easy to prove the following properties of $\leq_\sigma$:
\begin{enumerate}[(i)]
\setlength{\itemsep}{-5pt}
\item $\leq_\sigma$ defines a directed preorder on $M$.
\item $m \leq_\sigma m'$ and $m' \leq_\sigma m$ iff $m - m' \in \sigma_M^\bot$.
\item If $\tau \leq \sigma$, then $m \leq_\sigma m'$ implies $m \leq_\tau m'$.
\item If $\sigma$ is of maximal dimension in $N$, then $\leq_\sigma$ is a partial order.
\end{enumerate}

\begin{definition}
Let $\{E_m^\sigma\}_{m \in M}$ be a family of $k$-vector spaces. For each relation
$m \leq_\sigma m'$ let there be given a
vector space homomorphism $\chi^\sigma_{m, m'} : E^\sigma_m \longrightarrow E^\sigma_{m'}$
such that $\chi^\sigma_{m, m} = 1$ for all $m \in M$
and $\chi^\sigma_{m, m''} = \chi^\sigma_{m', m''} \circ \chi^\sigma_{m, m'}$
for each triple $m \leq_\sigma m' \leq_\sigma m''$.
We denote such data by  $\hat{E}^\sigma$ and call it a {\it $\sigma$-family}.
\end{definition}

\begin{lemma}
\label{stabrelevant}
Let $\hat{E}^\sigma$ be a $\sigma$-family. Let
$m, m'$ be given such that $m - m' \in \sigma_M^\bot$. Then the homomorphisms $\chi^\sigma_{m, m'}$
and $\chi^\sigma_{m', m}$ are isomorphisms.
\end{lemma}

\begin{proof}
The lemma follows from the fact that $\chi^\sigma_{m, m'} \circ \chi^\sigma_{m', m} =
\chi^\sigma_{m',m'} = 1$.
\end{proof}

\begin{definition}
Let $\hat{E}^\sigma$ and $\hat{F}^\sigma$ be two $\sigma$-families with vector space homomorphisms
$\chi^\sigma_{m, m'}$ and $\psi^\sigma_{m, m'}$, respectively. Then a {\it morphism}
$\hat{\phi}^\sigma$ from $\hat{E}^\sigma$ to $\hat{F}^\sigma$ is a set of
vector space homomorphisms $\{\phi_m^\sigma : E^\sigma_m \longrightarrow F^\sigma_{m}\}_{m \in M}$
such that $\phi_{m'}^\sigma \circ \chi^\sigma_{m, m'} =  \psi^\sigma_{m, m'} \circ \phi^\sigma_m$
for all $m, m' \in M$ with $m \leq_\sigma m'$.
\end{definition}

It is then clear that $\sigma$-families form a category. Moreover:

\begin{proposition}
\label{gradedisequivariant}
The following three categories are equivalent:
\begin{enumerate}[(i)]
\item\label{equivisgradedisfamilyi} equivariant quasicoherent sheaves over $U_\sigma$,
\item\label{equivisgradedisfamilyii} $M$-graded \mcoord{\sigma}-modules with morphisms of degree 0, and
\item\label{equivisgradedisfamilyiii} $\sigma$-families
\end{enumerate}
\end{proposition}

\begin{proof}
The Grothendieck functor $\tilde{\ }$ and the global section functor $\Gamma(U_\sigma, .)$
establish an equivalence of categories between the category of quasicoherent sheaves over
$U_\sigma$ and the category of \mcoord{\sigma}-modules. In the case that $U_\sigma$ is
smooth it was
shown by Kaneyama in \cite{Kan1} that the restrictions of these functors induce the
equivalence between (\ref{equivisgradedisfamilyi}) and (\ref{equivisgradedisfamilyii}).
Kaneyama's proof applies without change to the general situation.

We show the equivalence between (\ref{equivisgradedisfamilyii}) and
(\ref{equivisgradedisfamilyiii}).
From the discussion above it is clear that each $M$-graded \mcoord{\sigma}-module $E^\sigma$
gives rise to a $\sigma$-family $\hat{E}^\sigma$ via the family of vector spaces given by
the homogeneous components $E^\sigma_m$ and setting $\chi^\sigma_{m,m'}$ the map given
by multiplication with $\chi(m'-m)$. Moreover,
each homogeneous morphism $\phi$ by decomposition into homogeneous components
gives rise to a morphism $\hat{\phi}$ of $\sigma$-families.
On the other hand, given a $\sigma$-family $\hat{E}^\sigma$, we associated a graded
\mcoord{\sigma}-module to $\hat{E}^\sigma$ by setting
\begin{equation*}
E^\sigma := \bigoplus_{m \in M} E^\sigma_m
\end{equation*}
and for $m \in \sigma_M$ and a homogeneous element $e \in E^\sigma_{m'}$ we set
$\chi(m) \cdot e := \chi^\sigma_{m, m+m'}(e)$.
We then obtain the structure of a graded \mcoord{\sigma}-module by $k$-linear
continuation. Analogously, morphisms are composed by forming direct sums.
\end{proof}

\subsection{The Category of $\Delta$-Families}

Let $f : U_{\sigma'} \longrightarrow U_\sigma$ be an equivariant morphism of affine toric
varieties. This means that its restriction to the torus $T' \subset U_{\sigma'}$ is a morphism
of algebraic groups $T' \longrightarrow T \subset U_\sigma$ and that $f(t . x) = f(t) . f(x)$
for all $t \in T$ and $x \in U_{\sigma'}$.
Such an $f$ induces contravariantly a homomorphism of rings $f^* : \coord{\sigma}
\longrightarrow k[\sigma'_M]$ and a homomorphism $\check{f} : M \longrightarrow M'$
between the character groups of $T$ and $T'$, respectively. Both are compatible in the
sense that the restriction $\check{f} : \sigma_M \longrightarrow \sigma'_{M'}$ of $f$ is a
homomorphism of semigroups.
Then the pair $(f^*$, $\check{f})$ is a homomorphism of graded rings:
\begin{equation*}
(f^*, \check{f}) : (\coord{\sigma}, M) \longrightarrow (k[\sigma'_{M'}], M')
\end{equation*}
in the sense of Section \ref{gradings}.
Let \msh{E} be an equivariant quasicoherent sheaf on $U_\sigma$.
Then the sheaf $f^*\sh{E}$ is a $T'$-equivariant sheaf on $U_{\sigma'}$ and
the $k[\sigma'_{M'}]$-module $\Gamma(U_{\sigma'} f^*\sh{E})$ has an $M'$-grading,
where the natural isomorphism $\Gamma(U_{\sigma'} f^*\sh{E}) \cong \Gamma(U_\sigma, \sh{E})$
$\otimes_{\coord{\sigma}} k[\sigma'_{M'}]$ coincides with the graded scalar
extension in the sense of Definition \ref{gradscalext}.
The $M'$-graded $k[\sigma'_{M'}]$-module $\Gamma(U_{\sigma'} f^*\sh{E})$ thus corresponds
to a $\sigma'$-family.

\begin{definition}
Let $E$ be an $M$-graded \mcoord{\sigma}-module corresponding to a $\sigma$-family
$\hat{E}^\sigma$. Further let $f : U_{\sigma'} \longrightarrow U_\sigma$ be an equivariant
morphism of affine toric varieties. Then denote by $f^*\hat{E}^\sigma$ the $\sigma'$-family
obtained by the graded scalar extension of $E$ via the pair of homomorphisms $(f^*, \check{f}) :
(\coord{\sigma},M) \longrightarrow (k[\sigma'_{M'}], M')$.
\end{definition}

\begin{proposition}
Given an equivariant morphism $f : U_{\sigma'} \longrightarrow U_\sigma$ of affine
toric varieties. Then $f^*$ is up to natural equivalence a functor from the category of
$\sigma$-families to the category of $\sigma'$-families.
\end{proposition}

\begin{proof}
Note that $f^*$ in the category of sheaves is a functor only up to natural equivalence.
\end{proof}

We now paste together families of vector spaces associated to the cones of a fan:

\begin{definition}
\label{deltadef}
Let $\Delta$ be a fan. A collection $\{\hat{E}^\sigma\}_{\sigma \in \Delta}$ of
$\sigma$-families is called a {\it $\Delta$-family}, denoted $\hat{E}^\Delta$,
iff for each pair $\tau < \sigma$ with inclusions
$i^\tau_\sigma : U_\tau \hookrightarrow U_\sigma$
there exists an isomorphism of families $\eta_{\tau\sigma} : {i^\tau_\sigma}^*\hat{E}^\sigma
\overset{\cong}{\longrightarrow} \hat{E}^\tau$ such that for each triple
$\rho < \tau < \sigma$ there is an equality
$\eta_{\rho\sigma} = \eta_{\rho\tau} \circ {i^\rho_\tau}^*\eta_{\tau\sigma}$.
A {\it morphism} of $\Delta$-families is a collection of morphisms
$\{\hat{\phi}^\sigma : \hat{E}^\sigma \longrightarrow
\hat{E'}^\sigma\}_{\sigma \in \Delta}$ such that for all $\sigma, \tau$ and $\tau <
\sigma$ the following diagram commutes:
\begin{equation*}
\xymatrix{
{i^\tau_\sigma}^*\hat{E}^\sigma \ar[rr]^{{i^\tau_\sigma}^*\hat{\phi}^\sigma}
\ar[d]^{\eta_{\tau\sigma}} & & {i^\tau_\sigma}^*\hat{E'}^\sigma \ar[d]^{\eta'_{\tau\sigma}} \\
\hat{E}^\tau \ar[rr]^{\hat{\phi}^\tau} & & \hat{E'}^\tau\\
}
\end{equation*}
\end{definition}
With these definitions it is easy to see that the $\Delta$-families form a category.
Moreover:
\begin{theorem}
Let $\Delta$ be a fan. Then the category of $\Delta$-families
is equivalent to the category of quasicoherent equivariant sheaves over $X$.
\end{theorem}

\begin{proof}
Because the $\sigma$-families already have been shown to be equivalent to equivariant
quasicoherent sheaves on the invariant open subvarieties $U_\sigma$ of $X$, we have only to
show that Definition \ref{deltadef} encodes the usual data for gluing the sheaves
$\sh{E}_\sigma$ which correspond to the $\sigma$-families $\hat{E}^\sigma$ with respect
to the open cover $U_\sigma$, $\sigma \in \Delta$. This means that there exists a family
$\phi_{\sigma \sigma'}: \sh{E}_\sigma \vert_{U_{\sigma \cap \sigma'}} \overset{\cong}{\longrightarrow}
\sh{E}_{\sigma'} \vert_{U_{\sigma \cap \sigma'}}$ for all $\sigma, \sigma' \in \Delta$
such that $\phi_{\sigma \sigma} = \operatorname{id}$ and $\phi_{\sigma \sigma''} =
\phi_{\sigma' \sigma''} \circ \phi_{\sigma \sigma'}$ over each triple intersection
$U_{\sigma \cap \sigma' \cap \sigma''}$. Let $\sigma, \sigma' \in \Delta$ and set
$\tau := \sigma \cap \sigma'$ define
\begin{equation*}
\phi_{\sigma \sigma'} := \eta_{\tau \sigma'}^{-1} \circ \eta_{\tau \sigma}
\end{equation*}
Then the $\phi_{\sigma \sigma'}$ fulfill the cocycle condition up to natural equivalence
and the $\sh{E}_\sigma$ glue to a quasi-coherent sheaf \msh{E} on $X$. Furthermore it
is then straightforward to verify that $\hat{E}^\Delta \mapsto \sh{E}$ defines a functor
which induces an equivalence of categories.
\end{proof}

\paragraph{Coherence:}

The following finiteness conditions of a $\Delta$-family correspond to coherence.

\begin{definition}
\label{finiteness}
We say that a $\sigma$-family $\hat{E}^\sigma$ is {\it finite} iff
\begin{enumerate}[(i)]
\item\label{finii} all the vector spaces $E^\sigma_m$ are finite-dimensional,
\item\label{fini} for each chain $\dots <_\sigma m_{i - 1} <_\sigma
m_i <_\sigma \cdots$ of characters in $M$ there exists an $i_0 \in
\mathbb{Z}$ such that $E^\sigma_{m_i} = 0$ for $i < i_0$,
\item\label{finiii} there are only finitely many vector spaces $E_m$ such that the map
\begin{equation*}
\bigoplus_{m' <_\sigma m} E^\sigma_{m'} \longrightarrow E^\sigma_m
\end{equation*}
defined by summation of the $\chi^\sigma_{m',m}$, is not surjective.
\end{enumerate}
We say that a $\Delta$-family is {\it finite} if all of its $\sigma$-families are finite.
\end{definition}

\begin{proposition}
A quasicoherent equivariant sheaf is coherent iff its associated $\Delta$-family is finite.
\end{proposition}

\begin{proof}
Let $\hat{E}^\Delta$ be finite and let $\hat{E}^\sigma$ be any $\sigma$-family. Because the
homomorphisms $\chi^\sigma_{m, m'}$ are isomorphisms iff $m \leq_\sigma m'$ and
$m' \leq_\sigma m$, we assume without loss of generality that $\sigma$ is maximal
and thus $\leq_\sigma$ is a partial order.
Consider the set $P = \{m \in M / \sigma_M^\bot \text{ such that } E_m \neq 0\}$ and
let $P_{\min}$ be the set of minimal elements of $P$ with respect to $\leq_\sigma$.
Because of condition (\ref{finiii}), $P_{\min}$ must be finite, and by (\ref{fini}),
if the $\sigma$-family is nonempty, $P_{\min}$ is nonempty as well. Again by
(\ref{finiii}), there exist only finitely many cokernels
$\bigoplus_{m' <_\sigma m} E^\sigma_{m'} \longrightarrow E^\sigma_m$ which we denote
by $C_m$ and denote $Q := \{m \in M \setminus P \text{ such that } C_m \neq 0\}$.
The generators of the module $E^\sigma = \bigoplus_{m \in M} E^\sigma_m$ correspond
to the $k$-vector space basis of $\bigoplus_{m \in P_{min}} E^\sigma_m \oplus \bigoplus_{m \in Q}
C_m$ which by (\ref{finii}) is finite dimensional, hence $E^\sigma$ is finitely
generated. The converse follows straightforwardly by similar arguments.
\end{proof}

\begin{remark}
\label{stagnation}
Note that from Definition \ref{finiteness} it follows immediately that for every chain
$\dots <_\sigma
m_{i - 1} <_\sigma m_i <_\sigma \cdots$ the sequence $E^\sigma_{m_i}$ becomes
stationary, i.e. there exists
an $i_1 \geq i_0$ such that the $\chi^\sigma_{m_i, m_{i + 1}}$ become isomorphisms
for all $i \geq i_1$.
\end{remark}

\subsection{Torsion Free Equivariant Sheaves}

Let a finite $\Delta$-family $\hat{E}^\Delta$ be given.
Each of its $\sigma$-families $\hat{E}^\sigma$ is preordered, so there exists for each
$\sigma \in \Delta$ a direct limit $\dlim \hat{E}^\sigma$ which we denote by
$\mathbf{E}^\sigma$. Thus for each two vector spaces $E^\sigma_m, E^\sigma_{m'}$
in $\hat{E}^\sigma$ with $m \leq_\sigma m'$, there is a commutative diagram
\begin{equation*}
\xymatrix{
E^\sigma_m \ar[r]^{\chi^\sigma_{m,m'}} \ar[rd] & E^\sigma_{m'}
\ar[d]\\
& \mathbf{E}^\sigma
}
\end{equation*}

Recall that a coherent sheaf \msh{E} on a normal scheme $X$ is {\it torsion free} iff
there exists an injective homomorphism of $\sh{O}_X$-modules $\sh{E} \longrightarrow
k(X)^r$ for some $r \geq 0$, where $k(X)$ is the field of rational functions over $X$.
In our case \msh{E} being a torsion free equivariant sheaf
over a toric variety $X$, this implies that for all $\sigma \in \Delta$,
$\Gamma(U_\sigma, \sh{E})$ is a torsion free \mcoord{\sigma}-module. It is easy to see
that this is equivalent to the fact that $\chi(m) \cdot e \neq 0$ for all $m \in \sigma_M$
and all $M$-homogeneous elements $e$ of $\Gamma(U_\sigma, \sh{E})$. This in turn implies:

\begin{proposition}
Let \msh{E} be an equivariant coherent sheaf on $X$ and $\hat{E}^\Delta$
its $\Delta$-family.
\msh{E} is torsion free iff for all $\sigma \in \Delta$ the maps in the above diagram are injective.
\end{proposition}

If \msh{E} is torsion free, then the restriction map $\Gamma(U_\sigma, \sh{E})
\longrightarrow \Gamma(U_\tau, \sh{E})$ for any $\tau < \sigma$ is injective.
So there exists a natural inclusion of $\sigma$-families $\hat{E}^\sigma \hookrightarrow
(i^\tau_\sigma)^* \hat{E}^\sigma$ for any $\tau < \sigma$. Therefore the composition
\begin{equation}
\label{resinc}
\xymatrix{
\hat{E}^\sigma \ar@{^{(}->}[r] & (i^\tau_\sigma)^* \hat{E}^\sigma
\ar[r]_\cong^{\eta_{\tau\sigma}} & \hat{E}^\tau
}
\end{equation}
is injective. It is then easy to prove:

\begin{proposition}
Let \msh{E} be the coherent sheaf corresponding to the $\Delta$-family $\hat{E}^\Delta$. If
\msh{E} is torsion free then the homomorphisms $\hat{E}^\sigma \longrightarrow
\hat{E}^\tau$ as in (\ref{resinc}) are injections for any $\tau < \sigma$.
\end{proposition}

For the rest of this subsection we will assume that all sheaves in question are
torsion free.
The injection $\hat{E}^\sigma \hookrightarrow \hat{E}^\tau$ induces an injective map
$\tilde{\eta}_{\tau\sigma} : \mathbf{E}^\sigma \hookrightarrow
\mathbf{E}^\tau$. The system of vector spaces $\mathbf{E}^\sigma$ and homomorphisms
$\tilde{\eta}_{\tau\sigma}$ forms a directed partially ordered family, $\mathbf{E}^\Delta$,
with respect to the reversed partial order '$<$' among the cones of $\Delta$. For any
$\sigma_1, \sigma_2 \in \Delta$ we have diagrams
\begin{equation*}
\xymatrix{
\mathbf{E}^{\sigma_1} \ar@{^{(}->}[rd] \\
& \mathbf{E}^{\sigma_1 \cap \sigma_2} \ar@{^{(}->}[r] & \mathbf{E}^{0}\\
\mathbf{E}^{\sigma_2} \ar@{^{(}->}[ru]
}
\end{equation*}
where $0$ is the minimal cone in $\Delta$.
We obtain an identification of the direct limit $\dlim \mathbf{E}^\Delta$ with
$\mathbf{E}^0$. By Proposition \ref{maxnullfree} any $M$-graded module over the
ring \mcoord{\sigma} is free, i.e. isomorphic to $k[M]^r$, for some $r \geq 0$.
This means that every equivariant sheaf over the torus $T = U_0$ is free and
$r$ is the rank of this sheaf. Consequently,
the morphisms $\chi^0_{m,m'}$ are all isomorphisms and $\mathbf{E}^0$ is an $r$-dimensional
$k$-vector space. Recall that the restriction map $\Gamma(U_\sigma, \sh{E}) \hookrightarrow
\Gamma(U_\tau, \sh{E})$ can be identified with the canonical map $\Gamma(U_\sigma, \sh{E})
\hookrightarrow \Gamma(U_\sigma, \sh{E})_{\chi(m_\tau)}$ into the localization with respect
to some element $m_\tau$ in the interior of $\check{\sigma} \cap \tau^\bot$.
By Remark \ref{stagnation} we know that for each chain $\dots <_\sigma m + i \cdot m_\tau
<_\sigma m + (i + 1) \cdot m_\tau <_\sigma \cdots$ the sequence $E^\sigma_{m + i \cdot m_\tau}$
becomes stationary and because $m + (i + 1) \cdot m_\tau - (m + i \cdot m_\tau) = m_\tau
\in \tau^\bot_M$ the maps
$\chi^\tau_{m + i \cdot m_\tau, m + (i + 1) \cdot m_\tau}$ are isomorphisms for all
$i \in \mathbb{Z}$. Therefore we have

\begin{proposition}
Let $\tau < \sigma$ and $\hat{E}^\sigma$ a $\sigma$-family. Let $m_\tau$ be an
integral element of the interior of $\check{\sigma} \cap \tau^\bot$ such that $\tau_M =
\sigma_M + \mathbb{Z}_{\geq 0} (-m_\tau)$. For each $m \in M$ there is the chain $\dots
<_\sigma m + i \cdot m_\tau <_\sigma m + (i + 1) \cdot m_\tau <_\sigma \cdots$. Then there
is an $i_m \in \mathbb{N}$ such that $((i^\tau_\sigma)^*(\hat{E}^\sigma))_m
\cong E^\sigma_{m + i \cdot m_\tau}$ for all $i \geq i_m$.
\end{proposition}

Applying this to $0 < \sigma$ we get:

\begin{corollary}
\label{dimisrank}
The inclusions $\mathbf{E}^\sigma \hookrightarrow \mathbf{E}^0$ are isomorphisms
$\mathbf{E}^\sigma \cong \mathbf{E}^0$.
\end{corollary}

\subsection{The Category of Multifiltrations}

Now we collect all the properties derived for equivariant, coherent torsion free sheaves:

\begin{definition}
\label{torsdef}
Let $\Delta$ be a fan, $V$ a finite-dimensional $k$-vector space, and let for each
$\sigma \in \Delta$ a set of subvector spaces $\{E^\sigma_m\}_{m \in M}$ of $V$ be given.
We say that this system is a {\it family of multifiltrations of $V$} if:
\begin{enumerate}[(i)]
\item For $\sigma \in \Delta$ and $m \leq_\sigma m'$, $E^\sigma_m$ is contained in $E^\sigma_{m'}$.

\item $V = \bigcup_{m \in M} E^\sigma_m$ for any $\sigma \in \Delta$.

\item For each chain $\dots <_\sigma m_{i - 1} <_\sigma m_i <_\sigma
\cdots$ of characters in $M$ there exists an $i_0 \in \mathbb{Z}$ such that $E^\sigma_{m_i}
= 0$ for all $i \leq i_0$.

\item\label{classifour} There exist only finitely many vector spaces $E^\sigma_m$ which are not
contained in the union of all vector spaces $E^\sigma_{m'}$ with $m' <_\sigma m$.

\item\label{classifive} (compatibility condition) For each $\tau < \sigma$ with $\tau_M = \sigma_M +
\mathbb{Z}_{\geq 0} (-m_\tau)$ we consider with respect to the preorder $\leq_\sigma$
the ascending chains $m + i \cdot m_\tau$ for
$i \geq 0$. By condition (\ref{classifour}) and because $V$ is finite dimensional the
sequence of subvector spaces $E^\sigma_{m + i \cdot
m_\tau}$ necessarily becomes stationary for some $i^\tau_m\in \mathbb{Z}$. We require that $E^\tau_m
= E^\sigma_{m + i^\tau_m \cdot m_\tau}$ for all $m \in M$.
\end{enumerate}
\end{definition}

Families of multifiltrations are $\Delta$-families which are realized as subvector spaces of
the limit vector space $\mathbf{E}^0$. A morphism of families of multifiltrations
$\{E^\sigma_m\}_{\sigma \in \Delta, m \in M}$ and $\{F^\sigma_m\}_{\sigma \in \Delta, m \in M}$
then is equivalent to a homomorphism of vector spaces $\mathbf{E}^0 \longrightarrow \mathbf{F}^0$
which is compatible with these multifiltrations and so induces a morphism of $\Delta$-families.
This technical reformulation gives:

\begin{theorem}
The category of torsion free equivariant coherent sheaves is equivalent to the category of
families of multifiltrations of finite-dimensional vector spaces.
\end{theorem}

Now we briefly explain how the above classification is related to the characterization
of torsion free sheaves in \cite{Kly91}.
There it was stated that on a smooth complete toric variety $X$ the category
of equivariant torsion free sheaves is equivalent to the category of {\it multifiltrations}
for finite-dimensional vector spaces.
If $X$ is smooth, for any $r$-dimensional cone $\sigma$, its primitive vectors $n(\rho)$,
$\rho \in \sigma(1)$ form a part of a $\mathbb{Z}$-basis for $N$. Thus we can choose
a basis $\{m(\rho_i)\}_{\rho_i \in \sigma(1)}$ of $M / \sigma_M^\bot$
dual to the minimal submodule of $N$ which contains $\sigma \cap N$ and which is spanned by
the $n(\rho)$. With respect to the dual basis we write for an element $m \in M$ its
residual class $\overline{m} \in M / \sigma_M^\bot$ as $\overline{m} = \sum_{\rho_i \in
\sigma(1)} i_{\rho_i} \cdot m(\rho_i)$. So by identifying $\overline{m}$ with the tuple
$(i_{\rho_1}, \dots, i_{\rho_r})$, we identify  $M / \sigma_M^\bot$ with
$\mathbb{Z}^{\sigma(1)}$. Moreover, as subvector spaces of $\mathbf{E}^0$, we have that
$E^\sigma_m = E^\sigma_{m'}$ whenever $m - m' \in \sigma_M^\bot$. So we can write
$E^\sigma(i_{\rho_1},
\dots, i_{\rho_r})$ for $E^\sigma_m$.
Then a multifiltration in the sense of \cite{Kly91}
associated to a cone $\sigma$ is a set
$\{E^\sigma(i_{\rho_1}, \dots, i_{\rho_r})\}_{(i_{\rho_1}, \dots, i_{\rho_r})
\subset \mathbb{Z}^{\sigma(1)}}$
of subvector spaces of some vector space $\mathbf{E}^0$ which is parametrized by
indices $i_{\rho_k}$ for each $\rho_k \in
\sigma(1)$.
In this formulation the conditions for the filtrations are:
\begin{itemize}
\item $E^\sigma(i_{\rho_1}, \dots, i_{\rho_k}, \dots, i_{\rho_r})
\subset E^\sigma(i_{\rho_1}, \dots, i_{\rho_k} + 1, \dots, i_{\rho_r})$
for each $\rho_k \in \sigma(1)$.

\item $\bigcup_{i_{\rho_1}, \dots, i_{\rho_r}}
E^\sigma(i_{\rho_1}, \dots, i_{\rho_r}) = V$.

\item Let $\tau < \sigma$, then $\tau$ is spanned by the rays $\{\rho_{k_1}, \dots, \rho_{k_s}\}
= \tau(1) \subset \sigma(1)$. Assume that $(k_1, \dots, k_s) = (1, \dots, s)$, then
$E^\tau(i_{\rho_1}, \dots, i_{\rho_s}) = E^\sigma(i_{\rho_1}, \dots, i_{\rho_s}, \infty,
\dots, \infty)$ and analogously for $(k_1, \dots, k_s)$ different from $(1, \dots, s)$.
Here $\infty$ abbreviates the choice of a suitable $i^\tau_m \in \mathbb{Z}$ as in Definition
\ref{torsdef}, (\ref{classifive}).
\end{itemize}

These are the conditions of Klyachko (\cite{Kly91}).

\subsection{Reflexive Equivariant Sheaves}

If \msh{E} is a reflexive sheaf on a normal variety $X$, then $\Gamma(X, \sh{E}) =
\Gamma(X \setminus Y, \sh{E})$ if $Y$ is a closed subset of $X$ of codimension
at least two. If \msh{E} is an equivariant
reflexive sheaf on a toric variety $X = X_\Delta$, a natural choice for $Y$ is the union
of orbits in $X$ which have
codimension at least two, i.e. $\Gamma(X, \sh{E}) = \Gamma(X_{\Delta_1}, \sh{E})$,
where $\Delta_1 = \Delta(0) \cup \rays$. In particular, on an affine
toric variety we can write $\Gamma(U_\sigma, \sh{E}) = \Gamma\big(\bigcup_{\rho \in \sigma(1)}
U_\rho, \sh{E}\big)$. So, if we consider the $\Gamma(U_\rho, \sh{E})$ as $k$-subvector spaces
of $\Gamma(U_0, \sh{E})$, we have
\begin{equation*}
\Gamma(U_\sigma, \sh{E}) \cong \bigcap_{\rho \in \sigma(1)} \Gamma(U_\rho, \sh{E}).
\end{equation*}
This implies that $\Gamma(U_\sigma, \sh{E})_m \cong \bigcap_{\rho \in \sigma(1)}
\Gamma(U_\rho, \sh{E})_m$ for each graded component of degree $m$.
We can translate this in a natural way to the intersection of multifiltrations in
the limit vector space $\mathbf{E}^0$:
\begin{equation*}
E^\sigma_{m} = \bigcap_{\rho \in \sigma(1)} E^\rho_m.
\end{equation*}
Hence a reflexive sheaf is completely determined by the multifiltrations $E^\rho_m$ of $V$
for $\rho \in \rays$.
We know that $E^\rho_m = E^\rho_{m'}$ if $m - m' \in \rho^\bot_M$. Thus the multifiltrations
are determined by the stabilizer of the minimal orbit of $U_\rho$, whose group of characters
is $M / \rho^\bot_M$ which canonically can be identified with $\mathbb{Z}$  by using the
primitive lattice element $n(\rho)$.
So we can pass to a full filtration of $\mathbf{E}^0$ in the usual sense:
\begin{equation*}
0 \quad \ldots \subset E^\rho(i) \subset E^\rho(i + 1) \subset \dots \quad \mathbf{E}^0
\end{equation*}
and have the identity:
\begin{equation*}
E^\rho_m = E^\rho(\langle m , n(\rho) \rangle)
\end{equation*}
Using this identity we easily get back a $\Delta$-family from a given set of filtrations
$E^\rho(i)$ by setting
\begin{equation*}
E^\sigma_m := \bigcap_{\rho \in \sigma(1)} E^\rho(\langle m, n(\rho) \rangle ).
\end{equation*}
Now we have a classification for equivariant reflexive sheaves on toric varieties:

\begin{theorem}
The category of equivariant reflexive sheaves on a toric variety $X$ is equivalent to
the category of vector spaces with full filtrations associated to each ray in \mrays.
The morphisms in this category are vector space homomorphisms which are compatible with
the filtrations in the $\Delta$-family sense.
\end{theorem}

\paragraph{Locally free sheaves:}

\begin{proposition}
\label{affinefree}
Finitely generated projective $M$-graded \mcoord{\sigma}-modules are free.
\end{proposition}

\begin{proof}
The homogeneous ideal of \mcoord{\sigma} associated to the semigroup ideal $\sigma_M \setminus
\sigma_M^\bot$ is the unique maximal homogeneous ideal of \mcoord{\sigma}, so by
Corollary \ref{locgradprojfree}, projective $M$-graded \mcoord{\sigma}-modules are free.
\end{proof}

\begin{remark}
Note that  projective \mcoord{\sigma}-modules correspond to locally free sheaves over
$U_\sigma$. So by this proposition we obtain another proof for Proposition 2.1 i) in
\cite{Kly90} and Theorem 3.5 in \cite{Kan1}: equivariant vector bundles on affine toric
varieties are free.
\end{remark}

An equivariant bundle \msh{E} of rank $r$ on $U_\sigma$ corresponds to
the graded free \mcoord{\sigma}-module
\begin{equation*}
\Gamma(U_\sigma, \sh{E}) = \bigoplus_{i = 1}^r \coord{\sigma}(-m_i) =
\bigoplus_{i = 1}^r \chi(m_i) \cdot \coord{\sigma}
\end{equation*} 
i.e. it is a direct sum of graded free modules with respect to $r$ characters $m_i \in M$.
The latter equality comes from an equivariant embedding of $\Gamma(U_\sigma, \sh{E})$
into the quasi-coherent \mcoord{\sigma}-module $\Gamma(T, \sh{E}) \cong k[M]^r$.

Being reflexive, locally free sheaves are determined by 
a set of filtrations $\{E^\rho(i), i \in \mathbb{Z}\}_{\rho \in \rays}$ of some vector
space $E = \mathbf{E}^0$. By Proposition \ref{affinefree} any equivariant vector bundle
over an affine toric variety $U_\sigma$ is trivial and $\Gamma(U_\sigma, \sh{E})$
decomposes into a direct sum of free \mcoord{\sigma}--modules of rank 1:
\begin{equation*}
\Gamma(U_\sigma, \sh{E}) = \bigoplus_{i = 1}^r \coord{\sigma}(-m_i)
\end{equation*}
where $r$ is the rank of \msh{E}. So for any $m \in M$
\begin{equation*}
\Gamma(U_\sigma, \sh{E})_m = \bigcap_{\rho \in \sigma(1)} \Gamma(U_\rho, \sh{E})_m
\end{equation*}
Thus we arrive at Klyachko's compatibility condition for the filtrations $E^\rho(i)$:

\begin{theorem}[\cite{Kly90}]
\label{kly90}
The category of equivariant vector bundles over a toric variety $X_\Delta$ is equivalent to
the category of vector spaces with a family of filtrations $E^\rho(i)$ for each $\rho
\in \rays$ which satisfy the following compatibility condition:

Let $E$ be a vector space with a family of filtrations, then for any $\sigma
\in \Delta$ there exists a $T$-eigenspace decomposition $E = \bigoplus_{m \in M} E_m$
such that
\begin{equation*}
E^\rho(i) = \sum_{m, \langle m , n_\rho \rangle \leq i} E_m
\end{equation*}
\end{theorem}

Note that here we are considering {\it increasing} filtrations instead of decreasing
ones as in \cite{Kly90}.

\subsection{Equivariant Sheaves and Homogeneous Coordinate Rings}

As in Section \ref{homcoordring} let $C$ be the positive orthant of $\mathbb{R}^\rays$,
$\hat{\Delta}$ the subfan of $C$ defining $\hat{X} = X_{\hat{\Delta}}$ and let $S =
k[x_\rho \mid \rho \in \rays]$ be the homogeneous coordinate ring of X.
$S$ has a natural
\mweildivisors-grading, and we call a \mweildivisors-graded module over $S$ {\it fine-graded}.

\begin{proposition}[\cite{CoxBat}, Proposition 4.17]
\label{finegradedimpliesequivariant}\label{Zgrad}
\begin{enumerate}[(i)]
\setlength{\itemsep}{-5pt}
\item Let $E$ be a fine-graded $S$-module. Then the sheaf $\tilde{E}$ over $X$ has a natural
$T$-lineariza\-tion.
\item Let $X$ be an arbitrary toric variety. Then each equivariant sheaf $\sh{E}$ is
of the form $\tilde{E}$ for some fine-graded $S$-module $E$.
\end{enumerate}
\end{proposition}

Batyrev and Cox have proven the latter statement only for the simplicial case.

\begin{proof}
Consider the projection map $\pi : \hat{X} \longrightarrow X$ as in Section \ref{homcoordring}.
According to Corollary \ref{pullbacktilde}, $\sh{E}$ is isomorphic to
$\Gamma(\hat{X}, \pi^*\sh{E})\tilde{\ }$. Because $\pi$ is an equivariant map
of toric varieties, $\pi^*\sh{E}$ is a $(k^*)^\rays$-equivariant
sheaf on $\hat{X}$, hence by Proposition \ref{gradedisequivariant}
$\Gamma(\hat{X}, \pi^*\sh{E})$ possesses a
\mweildivisors-grading, and therefore is fine-graded.
\end{proof}

As we have seen in Section \ref{coxqcoh}, one obtains all possible reflexive sheaves of
rank 1 by shifting the $A_{n - 1}(X)$-degree of $S$. Analogously, by picking some $\underline{n}
\in \weildivisors$, which is equivalent to picking a $T$-invariant Weil divisor $D$ on $X$,
the module $\widetilde{S(\underline{n})}$ is naturally isomorphic to $\sh{O}_X(D)$ together
with its $T$-equivariant structure.

\begin{remark}
By \cite{ERR}, Theorem 1.1 we know that on a toric variety each sheaf has a resolution of
length $\leq \# \rays$ and that in the equivariant case this resolution can be chosen
in terms of $T$-linearized sheaves $\sh{O}_X(D)$ where $D$ is a $T$-invariant Weil divisor.
\end{remark}

A homogeneous morphism of fine-graded free modules
\begin{equation*}
\bigoplus_{j = 1}^m S(\underline{m}_j) \overset{A}{\longrightarrow}
\bigoplus_{i = 1}^n S(\underline{n}_i)
\end{equation*}
is given by an $n \times m$-matrix $A = (a_{ij})$ whose entries are monomials $a_{ij} =
\alpha_{ij}x^{\underline{m}_j - \underline{n}_i}$, where $\alpha_{ij} = 0$ whenever
$\underline{m}_j - \underline{n}_i \notin \mathbb{N}^\rays$.
The degree of the
monomials is completely determined by the grading. Hence, it is equivalent to
define such a morphism by a matrix $(\alpha_{ij})$ with $\alpha_{ij} \in k$.

A fine-graded $S$-module $E$ gives rise to an equivariant
quasicoherent sheaf over $k^\rays$ and is therefore equivalent to a family
of vector spaces $\hat{E}^C$.
One can easily see that $\hat{E}^C$ induces a $\hat{\Delta}$-family by restriction
to $\hat{X}$.
Over $U_{\hat{\sigma}}$ the module $E$ localizes to the \mcoord{\hat{\sigma}}-module
$E_{x^{\hat{\sigma}}}$.
Its zero component $E_{(x^{\hat{\sigma}})}$ is an $M$-graded $S_{(x^{\hat{\sigma}})}$-module.
By the inclusion $0 \longrightarrow M \overset{j}{\longrightarrow} \weildivisors$ from
diagram (\ref{maindiagram}) we can write
\begin{equation*}
E_{(x^{\hat{\sigma}})} = \bigoplus_{m \in M} (E_{x^{\hat{\sigma}}})_{j(m)}
\end{equation*}
Thus the set of homogeneous components $\{(E_{(x^{\hat{\sigma}})})_m\}_{m \in M}$
can be considered as a subset of the set of homogeneous components of
$E_{x^{\hat{\sigma}}}$. If we set $\chi^\sigma_{m, m'} := \chi^{\hat{\sigma}}_{j(m), j(m')}$
we obtain a $\sigma$-family as a subfamily of the $\hat{\sigma}$-family $\hat{E}^{\hat{\sigma}}$.

Thus we obtain a $\Delta$-family associated to a quasicoherent
equivariant sheaf over $X$ as a subfamily of the $\hat{\Delta}$-family
associated to the fine-graded module $E$.
We can refine this for torsion free equivariant sheaves as follows. Let $E$ be a
fine-graded $S$-module such that $\tilde{E}$ is a torsion free sheaf on $X$.
The torsion submodule $t(E)$ of $E$ is fine-graded as well and
the module $(E/t(E))$ gives rise to an equivariant sheaf $(E/t(E))\tilde{\ }$ over
$X$ isomorphic to $\tilde{E}$.
Thus if we consider torsion free sheaves of the form $\tilde{E}$ for some fine-graded
$S$-module $E$, we can assume that $E$ is torsion free and thus equivalent to a
family $E^{\hat{\Delta}}$ of multifiltrations for some vector space $V$ associated to the fan
$\hat{\Delta}$. By the above construction the $\Delta$-family obtained as subfamily of the
$\hat{\Delta}$-family defined by $E$
is a collection of families of subvector spaces of $V$. These families
are indeed multifiltrations of $V$ as well:

\begin{proposition}
\label{sublim}
$\dlim \hat{E}^{\hat{\sigma}} = \dlim \hat{E}^\sigma$ for all $\sigma \in \Delta$.
\end{proposition}

\begin{proof}
Because $\hat{E}^\sigma$ is a subfamily of $\hat{E}^{\hat{\sigma}}$ we know that
$\dlim \hat{E}^\sigma \subset \dlim \hat{E}^{\hat{\sigma}}$. For the other inclusion
it suffices to show that for each $\underline{n} \in \weildivisors$ there exists
$m \in M$ such that $\underline{n} <_{\hat{\sigma}} j(m)$. Denote by $\sup(\sigma)
\in \weildivisors$ the tuple whose $\rho$th entry is $1$
iff $\rho \in \sigma(1)$ and 0 else. Then the claim follows from the fact
that for an $m \in \operatorname{int}\check{\sigma} \cap
M$ and a suitable $n \in \mathbb{N}$ each component of $j(m) + n \cdot
\sup(\sigma)$ is bigger than zero. Then for each $\underline{n} \in \weildivisors$
there exists an $i \geq 0$ such that $\underline{n} <_{\hat{\sigma}} j(m) + (n + i) \cdot
\sup(\sigma)$.
It follows that $\underline{n} <_{\hat{\sigma}} j(m)$,
because with respect to the preorder $\leq_{\hat{\sigma}}$ it holds that $\underline{n'}
\leq_{\hat{\sigma}} \underline{n'} + k \cdot \sup(\sigma)$ and $\underline{n'} + k \cdot
\sup(\sigma) \leq_{\hat{\sigma}} \underline{n'}$ for all $\underline{n'} \in \weildivisors$
and all $k \in \mathbb{Z}$.
\end{proof}

\begin{example}
Consider a reflexive equivariant sheaf of rank 1 of the form $\widetilde{S(\underline{n})}$
for some $\underline{n} \in \weildivisors$.
Denote $\widetilde{S(\underline{n})}$ by $E$ and by $\underline{n}_\rho$ the component
of the tuple $\underline{n}$ corresponding to $\rho$.
$S(\underline{n})$ is equivalent to a set of filtrations for the 1-dimensional
vector space $k$:
\begin{equation*}
E^{\rho'}(i) =
\begin{cases}
0 & \text{ if } i < -\underline{n}_\rho \\
k & \text{ if } i \geq -\underline{n}_\rho
\end{cases}
\end{equation*}
for all $\rho' \in \Delta'(1)$. Each equivariant reflexive sheaf of rank 1 then is determined
by a function $\underline{n} : \rays \longrightarrow \mathbb{Z}$ (see also \cite{Kly90},
Example 2.3.1).
We can identify these filtrations with the filtrations for the corresponding reflexive sheaf on $X$:
\begin{equation*}
E^\rho(i) =
\begin{cases}
0 & \text{ if } i < -\underline{n}_\rho \\
k & \text{ if } i \geq -\underline{n}_\rho
\end{cases}
\end{equation*}
\end{example}

\begin{example}
Let $E$ be a fine-graded $S$-module such that $\sh{E}$ is a torsion free equivariant sheaf
of rank 1. As we have seen in the previous subsection, we can assume without
loss of generality that $E$ is a torsion free module of rank 1. Thus:

\begin{proposition}
$E$ is of the form $I(\underline{n})$, where $I$ is an ideal of $S$
generated by monomials.
\end{proposition}

By this proposition the classification of torsion free equivariant sheaves of rank
1 is reduced to the classification of monomial ideals. Such ideals have been investigated
(for example see \cite{monres1}) in terms of free resolutions. But the problem of finding
a minimal resolution for a general monomial ideal is still unsolved.
The one-to-one correspondence between invariant orbits of $\hat{X}$ and $X$ implies:

\begin{proposition}
\label{lci}
Let $\sigma \in \Delta$, then the ideal sheaf $\sh{I}_{V(\sigma)}$ is of the form
$\tilde{\mathfrak{I}}$ where $\mathfrak{I}$ is the \mweildivisors-homogeneous ideal
$\langle x_\rho \mid \rho \in \sigma(1) \rangle$ in $S$.
\end{proposition}
\end{example}

\section{Resolutions of Equivariant Vector Bundles of Rank Two on Smooth Toric Surfaces}
\label{resolutions}

Let $X$ be a smooth toric surface and let \msh{E} be a equivariant
vector bundle of rank 2 over $X$. Our goal is to construct a resolution for \msh{E}
under the assumption that \msh{E} is of sufficiently general type. The condition 'sufficiently
general' here is fulfilled by nearly all equivariant vector bundles of rank 2, so
that the knowledge of their resolutions allows a nearly exhaustive classification
of such vector bundles. The full classification problem for equivariant vector bundles
of rank 2 on toric surfaces will be the subject of a forthcoming paper.

Let us first show how equivariant vector bundles of rank 2 on a smooth toric surface
can be described in terms of filtrations. Let $\Delta$ be the fan describing a smooth
toric surface and let $n = \# \rays$
be the number of rays in $\Delta$. Then an equivariant vector bundle of rank 2 on
$X = X_\Delta$ is determined by a set of full filtrations $\dots E^\rho(i - 1) \subset
E^\rho(i) \subset E^\rho(i + 1) \subset \dots$ for $\rho \in \rays$, of a
2-dimensional vector space $\mathbf{E}^0$. A filtration $E^\rho$ is completely
described by two numbers $i^\rho_1 \leq i^\rho_2$ such that $0 = \dim E^\rho(i^\rho_1 - 1)
< \dim E^\rho(i^\rho_1)$ and $\dim E^\rho(i^\rho_2 - 1) < \dim E^\rho(i^\rho_2)$
and by one and the same vector space $E^\rho(i)$
for $i^\rho_1 \leq i < i^\rho_2$. Let us assume from now that $i^\rho_1 < i^\rho_2$, then
we can write a filtration as an ordered triple
$(i^\rho_1, i^\rho_2, E^\rho)$, where $i^\rho_1 \leq i^\rho_2$ and  $E^\rho \subset
E$ is a 1-dimensional subspace.
Because a reflexive sheaf on a smooth toric surface is always locally free, there are
no further
compatibility conditions (cf. Theorem \ref{kly90}) so that an equivariant rank 2 vector
bundle on $X$ generically is determined by an arbitrary collection of $2n$ numbers and $n$
1-dimensional subvector spaces of $\mathbf{E}^0$.

In \cite{Kan1}, Kaneyama classified equivariant vector bundles of rank 2 on
$\mathbb{P}_2$ using Euler sequences. Following this, we give a similar construction
idea of universal resolutions for
equivariant rank 2 bundles on arbitrary smooth toric surfaces, using
Euler
sequences for toric varieties as have earlier been constructed in \cite{Jaczewski}
and \cite{CoxBat}.

\begin{theorem}
\label{eulersequences}
Let $X$ be a smooth complete toric surface and let \msh{E} be an equivariant rank 2-vector
bundle over $X$ determined by filtrations $(i^\rho_1, i^\rho_2, E^\rho)$ as above. If
$i^\rho_1 < i^\rho_2$ for all $\rho \in  \rays$ and $E^\rho \neq V_{\rho'}$ whenever there
is a $\sigma \in \Delta(2)$ such that $\rho, \rho' \in \sigma(1)$, then there
exists, possibly after twisting by an equivariant line bundle, an Euler type equivariant short
exact sequence:
\begin{equation*}
0 \longrightarrow \sh{O}^{n - 2} \overset{A}{\longrightarrow} \bigoplus_{\rho \in \rays}
\sh{O}(i^\rho \cdot D_\rho) \longrightarrow \sh{E} \longrightarrow 0
\end{equation*}
where $n = \# \rays$, $i^\rho > 0$ and $A$ is an $n \times (n - 2)$-matrix of monomials.
\end{theorem}

\begin{proof}
Twisting an equivariant vector bundle \msh{E} by a $T$-linearized line bundle
$(S(\underline{n}))\tilde{\ }$ for some
$\underline{n} = (n_\rho) \in \weildivisors$ has the effect that the filtrations for \msh{E}
are shifted to $(i^\rho_1 + n_\rho, i^\rho_2 + n_\rho, E^\rho)$ for all $\rho$. Thus
we may assume that the filtrations are of the form $(-i^\rho, 0, E^\rho)$
for some positive integers $i^\rho$.

Consider now the quotient presentation $\hat{X} \longrightarrow X$.
Let $\{n_\rho\}_{\rho \in \rays}$ be the standard basis of the lattice $\hat{N}$. Then
a morphism of fine-graded free $S$-modules
\begin{equation*}
0 \longrightarrow S^{n-2} \overset{A}{\longrightarrow} \bigoplus_{\rho \in \rays} S(i^\rho n_\rho)
\end{equation*}
is defined by a matrix $A$. $A$ has monomial entries:
\begin{equation*}
A = (\alpha_{\rho,i} \cdot x_\rho^{i^\rho})
\end{equation*}
where $i$ runs from $1$ to $n - 2$. We require that $A$ has rank $n - 2$ as a homomorphism
of free $S$-modules.
After applying the sheafification functor $\tilde{\ }$
to this sequence we obtain a short exact sequence of sheaves
\begin{equation*}
0 \longrightarrow \sh{O}^{n - 2} \overset{A}{\longrightarrow} \bigoplus_{\rho \in \rays}
\sh{O}(i^\rho \cdot D_\rho) \longrightarrow \sh{Q} \longrightarrow 0
\end{equation*}
where by abuse of notion we write $A$ instead of $\tilde{A}$.

The matrix $A$ defines an injective morphism of coherent sheaves, but the image of
$\sh{O}^{n - 2}$ is not
necessarily a subbundle of $\bigoplus_{\rho \in \rays} \sh{O}(i^\rho \cdot D_\rho)$.
$A$ defines an inclusion of vector bundles if and only
if the rank of $A(x)$ equals $n - 2$ at all points
$x \in X$. This in turn means that $A$ is an inclusion of vector bundles
if and only if there exists a $k > 0$ such that $B^k \subset \fitt{A}$, where $B$
is the irrelevant ideal (see Section \ref{homcoordring}). If this is
the case, then the cokernel \msh{Q} is a vector bundle as well.
For $\phi, \psi \in \rays$, $\phi \neq \psi$ we introduce the following notation:
\begin{equation*}
x^{\phi\psi} :=\prod_{\rho \in \rays \setminus \{\phi,\psi\}}x_\rho^{i^\rho}.
\end{equation*}
Note that if $\sigma \in \Delta(2)$, and $\{\phi, \psi\} = \sigma(1)$ and $i^\rho =1$, for
$\rho \in \rays \setminus \{\phi,\psi\}$,
then $x^{\phi\psi} = x^{\hat{\sigma}}$ as defined in Section \ref{homcoordring}.
Let $A^{\phi\psi}$ be the $(n-2) \times (n-2)$-minor of $A$ which does not contain
the rows corresponding to $\phi$ and $\psi$ and let
\begin{equation*}
A' := (\alpha_{\rho,i})
\end{equation*}
be the matrix of coefficients of $A$ and $(A')^{\phi\psi}$ the corresponding minor.
The second Fitting ideal $\fitt{A}$ of $A$ is generated by the determinants of
all the $A^{\phi \psi}$:
\begin{equation*}
\fitt{A} = \langle \det A^{\phi\psi} \rangle =
\langle \det (A')^{\phi\psi} \cdot x^{\phi\psi} \rangle
\end{equation*}
Thus $\fitt{A}$ is a monomial ideal generated by those $x^{\phi\psi}$ for which $\det
(A')^{\phi\psi} \neq 0$. The condition $B^k \subset \fitt{A}$ for some $k > 0$ then is
equivalent to $\det (A')^{\phi\psi} \neq 0$ whenever $\phi$ and $\psi$ are contained in the
same 2-dimensional cone.

Denote $\mathbf{F}^0$ the $(n - 2)$-dimensional filtered $k$-vector space associated to the
vector bundle $\sh{O}^{n-2}$, and $\mathbf{G}^0$ the $n$-dimensional $k$-vector space
associated to $\bigoplus_{\rho \in \rays} \sh{O}(i^\rho \cdot D_\rho)$ which we identify with
$k^{\rays}$. Their filtrations are:
\begin{equation*}
F^\rho(i) =
\begin{cases}
0 & \text{\quad for } i < 0 \\
F & \text{\quad otherwise}
\end{cases}
\text{\quad and }
G^\rho(i) = 
\begin{cases}
0 & \text{\quad for } i < -i^\rho \\
k \cdot e_\rho & \text{\quad for } -i^\rho \leq i < 0 \\
\mathbf{G}^0 & \text{\quad otherwise}
\end{cases}
\end{equation*}
where the $e_\rho$ denote the standard basis of $k^\rays$. The matrix $A$ induces a
linear map of the vector spaces $\mathbf{F}^0$ and $\mathbf{G}^0$ which naturally
can be identified with the matrix $A'$.
We can define filtrations for the quotient vector space
$\mathbf{E}^0 := \mathbf{G}^0 / \mathbf{F}^0$ simply by taking the quotient filtrations
\begin{equation*}
E^\rho(i) = G^\rho(i) / F^\rho(i)
\end{equation*}
with respect to $A'$. These filtrations are of the form
\begin{equation*}
E^\rho(i) = (-i^\rho, 0, k \cdot \overline{e}_\rho)
\end{equation*}
where $\overline{e}_\rho$ is the image of $e_\rho$ in $\mathbf{E}^0$. If we assume that
$B^k \subset
\fitt{A}$ for some $k > 0$ then these filtrations become in a natural way the filtrations
associated to the cokernel \msh{E} of $A$.
On the other hand, if we define a homomorphism from $\mathbf{G}^0$ to some 2-dimensional
$k$-vector space $\mathbf{E}^0$ by fixing the
images $\overline{e}_\rho \neq 0$ of the basis vectors of $\mathbf{G}^0$, we immediately
obtain a homomorphism
of filtered vector spaces with kernel a filtered vector space
$\mathbf{F}^0 \overset{A'}{\hookrightarrow} \mathbf{G}^0$.
The corresponding matrix $A$ with monomial entries then defines a sheaf homomorphism
$0 \longrightarrow \sh{O}^{n - 2} \overset{A}{\longrightarrow} \bigoplus_{\rho \in \rays}
\sh{O}(i^\rho \cdot D_\rho)$. As we have seen before, the cokernel of $A$ is locally free
if and only if $\det(A)^{\phi\psi} \neq 0$ for all $\{\phi, \psi\} = \sigma(1)$ and
$\sigma \in \Delta(2)$.
Now it is a lemma from linear algebra that $\overline{e}_\phi$ and $\overline{e}_\psi$ are
linearly independent if and only if $\det (A')^{\phi\psi} \neq 0$.
\end{proof}

\end{document}